\documentclass[12pt]{article}
\usepackage{group21}
\title{Deformation quantization of compact K\"ahler manifolds
via Berezin-Toeplitz operators}
\author{Martin Schlichenmaier}
\address{Fakult\"at f\"ur Mathematik und Informatik,
        Universit\"at Mannheim,\\
        D-68131 Mannheim, Germany; schlichenmaier@math.uni-mannheim.de}
\usepackage{amsmath}
\usepackage{amsfonts}
\usepackage{amsthm}
\def\P{\mathbb P }

\def\Z{\mathbb Z }
\def\C{\mathbb C }

\def\w{\omega}
\def\i{{\,\mathrm i\,}}
\def\im{\text{Im\kern1.0pt }}
\def\re{\text{Re\kern1.0pt }}

\def\hh{\hat h}
\def\Tmp#1{T^{(m)}_{#1}}

\def\Pmm{{\mathcal{P}}(M)}

\def\Cim{C^{\infty}(M)}
\def\ghm{\Gamma_{hol}(M,L^{m})}
\def\gh{\Gamma_{hol}(M,L)}

\def\gul{\Gamma_{\infty}(M,L)}

\def\Lp{{\mathrm{L}}^2(M,L)}

\def\Qm#1{{Q}_{#1}^{(m)}}
\def\Qmo{{Q}^{(m)}}
\def\Pm#1{{P}_{#1}^{(m)}}

\def\Tfm{T_f^{(m)}}
\def\Tgm{T_g^{(m)}}
\def\Tfgm{T_{\{f,g\}}^{(m)}}

\def\zb{\overline{z}}

\def\d{\partial}
\def\db{\overline{\partial}}

\def\Pfz#1{\frac {\partial #1}{\partial z}}
\def\Pfzb#1{\frac {\partial#1}{\partial\overline{z}}}
\def\End{\text{End}}

\def\BT{Berezin-Toeplitz}

\def\CMP{Commun. Math. Phys. }
\def\JMP{Jour.  Math. Phys. }
\def\Izv{Math. USSR Izv. }
\newtheorem{theorem}{Theorem}
\theoremstyle{remark}
\newtheorem{example}{Example}
\newtheorem{remark}{Remark}

\begin{document}
\hspace*{\fill} Mannheimer Manuskripte 218

\hspace*{\fill} q-alg/9611022
\vspace{3.5truecm}
\begin{center}
\Large\bf 
Deformation quantization of 
\end{center}
\begin{center}
\Large\bf 
compact K\"ahler manifolds
\end{center}
\begin{center}
\Large\bf 
via Berezin-Toeplitz operators
\end{center}
\normalsize \rm
\begin{center}
{Martin Schlichenmaier}
\end{center}
\vspace{-0.8cm}
\begin{center}
{November 96, Manuskripte Nr. 218}
\end{center}
\vspace{5cm}
{Talk at the XXI Int. Coll. on Group Theoretical Methods in Physics, 
15-20 July, Goslar, Germany}
\maketitle
\begin{abstract}
This talk  reports  on results on the deformation quantization
(star products)  and on approximative operator representations for
quantizable compact K\"ahler manifolds obtained
via \BT\ operators. After choosing a holomorphic quantum line bundle
the \BT\ operator associated to a differentiable function on the 
manifold is the operator defined  by multiplying  global
holomorphic sections of the line bundle with this  function and
projecting the  differentiable section back to
the subspace of holomorphic sections.
The results were obtained in (respectively based on)
 joint work with M. Bordemann and E. Meinrenken.
\end{abstract}

\section{The set-up}
Let $\ (M,\w)\ $ be a compact (complex) K\"ahler manifold.
It should be considered as  phase space manifold $M$
with symplectic form given by the K\"ahler form $\w$.
Note that by the assumed compactness  the ``free phase space''
is not included here. Nevertheless there exist important
examples of such compact situations. E.g.~they could  
appear as
phase space manifolds for certain constrained systems and they could
appear after
dividing out certain symmetry group action in a
noncompact situation.
They also play a role in the quantization of
2-dimensional conformal field theory. Here the  space
which has to be quantized is the compactified moduli space of
stable holomorphic vector bundles (with additional structures)
on a Riemann surface. 

Using the K\"ahler form 
one can assign to every  differentiable function\footnote{
Differentiable always means arbitrary often
differentiable.} $f$ its Hamiltonian vector field $X_f$ and 
to every pair of functions $f$ and $g$ its
Poisson bracket:
\begin{equation}\label{lPoi}
\w(X_f,\cdot)=df(\cdot),\qquad  
\{\,f,g\,\}:=\w(X_f,X_g)\ .
\end{equation}
With the Poisson bracket the algebra of differentiable functions
$\Cim$ 
becomes a Poisson algebra $\Pmm$.
Assume $(M,\w)$ to be quantizable. This says there exists an
associated quantum line bundle $(L,h,\nabla)$.
Here $L$ is a holomorphic line bundle, $h$ a Hermitian metric on $L$ and 
$\nabla$ a connection compatible with the complex structure and 
the Hermitian metric.
With respect to a 
local complex coordinate patch and a 
local holomorphic frame of the bundle the metric $h$  is given by
a function $\hh$ and the connection as  ($d=\d+\db$)
\begin{equation}\label {lcon}
\nabla_{|}=\d+\d\log\hh+\db\  .
\end{equation}
The quantization condition  says that the 
curvature $F$ of the line bundle and the K\"ahler form $\w$ of the manifold
are related  as 
\begin{equation}\label{lquant}
F(X,Y):=\nabla_X\nabla_Y-\nabla_Y\nabla_X-\nabla_{[X,Y]}
=
-\i\w(X,Y)\ .
\end{equation}
In terms of the metric $h$ this can be written as $\ \i\db\d\log\hh=\w\ $.
The quantization condition  implies that $L$ is a positive (or ample)
line bundle. By the Kodaira embedding theorem  a certain
tensor power of $L$ is very ample, i.e. its  global
sections can be used to embed the phase space manifold $M$ 
into projective
space (with dimension given by the Riemann-Roch formula).
Hence, quantizable compact K\"ahler manifolds are projective
algebraic manifolds (and vice versa, see below).
\begin{example}
The Riemann sphere, the complex projective line,
 $\P(\C)=\C\cup \{\infty\}\cong S^2$. With respect to the quasi-global 
coordinate $z$  the form $\w$ can be given as
\begin{equation}
\w=\frac {\i}{(1+z\zb)^2}dz\wedge d\zb\ .
\end{equation}
The quantum line bundle  $L$ is the hyperplane  bundle.
For the Poisson bracket one obtains
\begin{equation}
\{f,g\}=\i(1+z\zb)^2\left(\Pfzb f\cdot\Pfz g-\Pfz f\Pfzb g\right)\ .
\end{equation}
\end{example}
\begin{example}
The (complex-) one dimensional torus $M$.
Up to isomorphy it can be given as
$M\cong \C/\Gamma_\tau$, where $\ \Gamma_\tau:=\{n+m\tau\mid n,m\in\Z\}$
is a lattice with $\ \im \tau>0$.
As K\"ahler form we take
\begin{equation}
\w=\frac {\i\pi}{\im \tau}dz\wedge d\zb\ ,
\end{equation}
with respect to the coordinate $z$ on the
covering space $\C$.
The corresponding quantum line bundle is the theta line bundle
of degree 1, i.e. the bundle whose global sections are
the scalar multiples of the Riemann theta function.
\end{example}
\begin{example}
The complex projective space $\P^n(\C)$ with
the Fubini-Study fundamental form as
K\"ahler form.
The associated quantum bundle is the hyperplane bundle.
By restricting the form and the bundle on
complex submanifolds (which are automatically algebraic)
all projective manifolds are quantizable.
As remarked  above every quantizable compact K\"ahler  manifold is 
projective. Nevertheless the embedding into 
projective space does not necessarily respect  
the  K\"ahler form.
\end{example}
For the space of  (arbitrary) global sections of $L$ a scalar product
is defined by
\begin{equation}\label{lsca}
\langle s_1,s_2\rangle=\int_M h(s_1,s_2)(x)\Omega(x),\qquad
\Omega=\frac 1{n!}
\underbrace {\w\wedge \w\cdots\wedge \w}_{\text{ $\dim_{\Bbb C} M$ times}}
\  .
\end{equation}
\section{Approximation results  for Toeplitz operators}
Let $\gul$ be the  (infinite-dimensional) space of global  
differentiable sections of the line bundle $L$,
$\Lp$ its L${}^2$-completion with respect to the scalar product
(\ref{lsca}), and
$\gh$ the finite-dimensional subspace
of global holomorphic sections.
Denote by $\Pi:\Lp\to\gh$ the projection.
The Toeplitz operator associated to the function $f\in\Pmm$ is
defined as the composition
\begin{equation}\label{lToe}
T_f:\gh\xrightarrow{M_f}\gul\xrightarrow{\Pi}\gh,\quad
s\mapsto f\cdot s\mapsto \Pi(f\cdot s)=:T_f(s)\ .
\end{equation}
Here $M_f$ denotes  the multiplication of the sections by
the function $f$.

Due to the  finite-dimensionality of $\gh$
information gets lost.
To recover this information we have to do everything
for each tensor power
$\ L^m=L^{\otimes m}$ of $L$.
This gives us
$\ (L^m,h^{(m)},\nabla^{(m)})$ and the  Toeplitz operators
$\Tmp f$.
The assignment
\begin{equation}\label{lBT}
T^{(m)}:\Pmm\to\End(\ghm),\qquad f\mapsto \Tfm
\end{equation}
is called {\em \BT\ quantization} (of level $m$).
The theorems below will justify the use of  the term ``quantization''. 
\begin{remark}
If one switches to the notation $\hbar=1/m$ then
the limit $m\to\infty$ corresponds to $\hbar\to 0$ as 
``classical limit''.
\end{remark}
\begin{remark}
By a result of Tuynman 
\cite{rTuyQ}
suitable  reinterpreted (see \cite{rBHSS} for a coordinate independent proof)
one obtains the relation 
$\  \Qm f=\i\Tmp {f-\frac {1}{2m}\Delta f}\ $,
where $\Qmo$ is the more well-known operator
of
geometric quantization with respect to the
prequantum operator
$\ \Pm f= -\nabla_{X_f^{(m)}}^{(m)}+\i f\cdot id\ $ and 
with K\"ahler polarization.
\end{remark}
\begin{remark}
In the example of 2-dimensional conformal field theory 
mentioned at the beginning, there exists a natural quantum 
line bundle $L$.
The spaces $\ghm$ are the {\em Verlinde spaces}. 
\end{remark}

Assume $L$ to be already very ample. By a 
rescaling of the K\"ahler form this always can be achieved.
The following two theorems are proved in 
a joint paper with M. Bordemann and E. Meinrenken 
\cite{rBMS}.
\begin{theorem}
For every  $\ f\in \Cim\ $ we have $\ ||\Tfm||\le||f||_\infty $ 
and
\begin{equation}\label{lThma}
\lim_{m\to\infty}||\Tfm||=||f||_{\infty}\ .
\end{equation}
Here $||f||_\infty$ is the sup-norm of $\ f\ $ on $M$ and
$||\Tfm||$ is the operator norm with respect to the 
scalar product (\ref{lsca}) on $\ghm$.
\end{theorem}
\begin{theorem}
For every  $\ f,g\in \Cim\ $ we have 
\begin{equation}\label{lThmb}
||m\i[\Tfm,\Tgm]-\Tfgm||\quad=\quad O(\frac 1m)\quad
\text{as}\quad m\to\infty
\ .
\end{equation}
\end{theorem}
By interpreting $1/m$ as $\hbar$ we see that up to order $\hbar^2$ the
quantization prescription
{\it ``Replace  the classical Poisson bracket
by the (rescaled) commutator of  the operators''}
is fulfilled. We obtain an approximative 
operator representation.
\section{Deformation quantization (star products)}
By a refinement of the methods 
for  the proofs of the above theorems it is possible to
construct a star product, resp. a  deformation
quantization for $\Pmm$.
Let me recall the definition of a star product.
Let $\mathcal{A}=\Cim[[\hbar]]$ be the algebra of formal power
 series in the
variable $\hbar$ over the algebra $\Cim$. A product $\star$
 on $\mathcal{A}$ is 
called a (formal) star product if it is an
associative $\C[[\hbar]]$-linear product such that
\begin{enumerate}
\item
$\mathcal{A}/\hbar\mathcal{A}\cong\Cim$, i.e.\quad $f\star g
 \bmod \hbar=f\cdot g$,
\item
$\dfrac 1\hbar(f\star g-g\star f)\bmod \hbar = -\i \{f,g\}$.
\end{enumerate}

We can always write 
\begin{equation}\label{lcif}
 f\star g=\sum\limits_{j=0}^\infty C_j(f,g)\hbar^j\quad\text{with}
\quad C_j(f,g)\in\Cim\ .
\end{equation}
The conditions 1. and 2.  can 
be reformulated as 
\begin{equation}\label{lcifa}
C_0(f,g)=f\cdot g,\qquad
C_1(f,g)-C_1(g,f)=-\i \{f,g\}\ .
\end{equation}
\begin{theorem}
There exists a unique (formal) star product on $\Cim$
\begin{equation} \label{stareq}
f \star g:=\sum_{j=0}^\infty \hbar^j C_j(f,g),\quad C_j(f,g)\in
C^\infty(M),
\end{equation}
in such a way that for  $f,g\in\Cim$ and for every $N\in\Bbb N$  we have
\begin{equation}
||T_{f}^{(m)}T_{g}^{(m)}-\sum_{0\le j<N}\left(\frac 1m\right)^j
T_{C_j(f,g)}^{(m)}||=K_N(f,g) \left(\frac 1m\right)^N
\end{equation}
for  $m\to\infty$, with suitable constants $K_N(f,g)$.
\end{theorem}
Again this is joint work with M. Bordemann and E. Meinrenken.
The proof can be found in \cite{rSchlH}.
Also, it will appear  in a forthcoming paper.
\section{Some additional remarks}
Time  and space does not permit me to give even  a rough sketch of the
proof.
Let me only point out the following steps. We
use the embedding of $M$
via the  sections of $L$ (which is now assumed to be very ample)
into projective space.
For the embedded situation we have an associated  Toeplitz structure,
the generalized Toeplitz operators
(of arbitrary orders) and their symbol calculus as
introduced by Boutet de  Monvel and Guillemin \cite{rBoGu}.
The above defined $\Tfm$ can be considered as the ``modes''
(with respect to an $S^1$-action) of a global Toeplitz operator
 $\tilde T_f$.
The symbol calculus gives bounds for their norms.

With the help of the above  theorems it is possible to show 
\cite{rBMS}
the
conjecture of Bordemann, Hoppe, Schaller  and Schlichenmaier
\cite{rBHSS} that for quantizable K\"ahler manifolds the Poisson
 algebra
is an $\mathfrak u(d(N))$, $N\to\infty$ quasi-limit with  a 
strictly increasing integer valued function 
$d(N)$  given by the Riemann-Roch formula.
For the  definition of quasi-limits 
I refer to the above mentioned article.

References to related works
can be found in \cite{rSchlB}. Here I only want to quote some names:
F.A.~Berezin, B.V.~Fedosov, J.H.~Rawnsley, M.~Cahen, S.~Gutt,
S.~Klimek, A.~Lesniewski, L.A.~Coburn,...
Note also that the Berezin-Toeplitz quantization           
can be considered in the frame-work of ``prime quantization''
as introduced  by S.T.~Ali and H.D.~Doebner.

The relation of \BT\ 
quantization  to Berezin's coherent states \cite{rBer} and symbols 
(in the coordinate independent formulation due to Rawnsley)
is studied in \cite{rSchlH}. In particular, it is
shown that for compact K\"ahler manifolds the Toeplitz map
$\ f\mapsto T_f\ $ and Berezin's covariant symbol map
$\ B\mapsto \sigma(B)\ $ are adjoint maps with respect to
the Hilbert-Schmidt scalar product for the operators  and 
the Liouville  measure on $M$ which is modified by the
 $\epsilon$-function of Rawnsley for the functions.

\end{document}